\documentclass[10pt, psamsfonts]{amsart}

\usepackage[a4paper,margin=4cm]{geometry}
\usepackage{amsfonts,amssymb,amscd}
\usepackage{graphicx}
\usepackage{subfigure}
\newtheorem{theorem}{Theorem}[section]

\newtheorem{lemma}[theorem]{Lemma}

\newtheorem{proposition}[theorem]{Proposition}
\newtheorem{example}{Example}

\theoremstyle{remark}

\newtheorem{remark}[theorem]{Remark}

\renewcommand{\leq}{\leqslant}
\renewcommand{\geq}{\geqslant}
\renewcommand{\le}{\leqslant}

\newcommand{\ptl}{\partial}

\newcommand{\ga}{\gamma}

\newcommand{\Om}{\Omega}
\newcommand{\sg}{\sigma}
\newcommand{\escpr}[1]{\left< #1\right>}

\newcommand{\rr}{\mathbb{R}}
\newcommand{\sph}{\mathbb{S}}
\newcommand{\nn}{\mathbb{N}}


\numberwithin{equation}{section}

\begin{document}

\title[Stable regions in rotationally symmetric tori]{Stable and
isoperimetric regions in rotationally symmetric tori with decreasing Gauss
curvature}

\author[A. Ca\~nete]{Antonio Ca\~nete}
\address{Departamento de Matem\'aticas \\ Campus de Ponferrada \\
Universidad de Le\'on \\ E-24400 Le\'on (Espa\~na)} 
\email{antonioc@unileon.es}

\date{April 6th, 2006}
\thanks{Work partially supported by MCyT-Feder research
project MTM2004-01387} 
\keywords{Stability, isoperimetric problem}
\subjclass[2000]{49Q20, 49Q10}


\begin{abstract}
\noindent In this work we classify the stable regions
(second order minima of perimeter under an area constraint)
in tori of revolution with piecewise continuous decreasing
Gauss curvature from the longest parallel
and with a horizontal symmetry.
Some applications to isoperimetric problems are also given.
\end{abstract}

\maketitle

\thispagestyle{empty}
\section*{Introduction}
\label{sec:introduction}

In a Riemannian surface $M$ we may consider the isoperimetric
problem consisting of finding the least perimeter sets in $M$
enclosing a given area $a_0$, with $a_0\leq area(M)$. If such a
set exists, then it is called an \emph{isoperimetric region}. In
the last years this problem has been of great interest, but only
for certain surfaces isoperimetric regions have been completely
classified (see \cite{bc},~\cite{p},~\cite{t},
~\cite{hhm1},~\cite{hhm2},~\cite{ritore}). A very important
related concept is the one of stability: a \emph{stable region} is a second
order minimum of perimeter under any variation preserving the
area enclosed. Variational formulae for perimeter and area imply that the
boundary of a stable region is composed of curves with the
\emph{same} constant geodesic curvature. Since any
isoperimetric region is stable, the characterization of the
stable regions in a surface is an interesting question in this setting.
Moreover, from a physical point of view, stable regions are
more realistic models since they are just local minima of
perimeter, instead of global as it is the case for an
isoperimetric region.

In this paper we deal with these problems in rotationally
symmetric tori with decreasing Gauss curvature from the
longest parallel and with a horizontal symmetry. This is a
large family of surfaces, including the standard tori
obtained by rotating a circle in $\rr^3$ with respect to
a line contained in the same plane as the circle,
and at a certain distance, and certain round spheres
to which hyperbolic annuli have been added. We obtain all possible stable
regions that may appear, checking as well if they occur as
isoperimetric.

The study of the above questions in surfaces of revolution with
decreasing Gauss curvature has been treated in different works.
Benjamini and Cao~\cite{bc} proved that in planes with total
positive curvature less than or equal to $2\pi$ and with
Gauss curvature decreasing from a pole, the isoperimetric solutions are
geodesic disks centered at the pole.
Later, Morgan, Hutchings and Howards~\cite{hhm2}
solved the problem in the general case of decreasing Gauss
curvature for planes, and also for real projective planes, annuli
with an end of finite area and certain spheres, obtaining that the
solutions are geodesic disks or annuli.

The proofs of the results in \cite{bc},~\cite{p},~\cite{t} and
~\cite{hhm2} involve different isoperimetric inequalities, which
are relations between the area and the perimeter of a set.
Ritor{\'e} uses in \cite{ritore} another approach, after Schmidt~\cite{s},
studying the closed embedded curves with constant geodesic curvature.
Since the boundary of any isoperimetric set consists of curves of this kind,
this technique allows to solve, after classifying stable regions,
the isoperimetric problem in planes of revolution, spheres with an
equatorial symmetry and projective planes, with decreasing and
also increasing curvature, and in certain annuli with decreasing
curvature.

In this paper we shall follow this second approach in order to classify the stable
regions in rotationally symmetric tori with decreasing Gauss
curvature from the longest parallel and with a horizontal
symmetry, and we obtain in our main classification Theorem~\ref{th:main}
disks bounded by constant geodesic curvature
curves (symmetric with respect to such parallel),
annuli bounded by two circles of revolution (symmetric or
nonsymmetric with respect to the shortest parallel), unions of
vertical annuli (each one bounded by two vertical geodesics), domains whose
boundary is an unduloid type curve and a circle of revolution, and
regions consisting of the union of a disk and a symmetric
annulus. We also check that all these regions appear as stable
ones in certain given surfaces.

Furthermore, we apply this classification to study the isoperimetric
problem described above. One of the most interesting consequences is that
an unduloid type curve may be part of the boundary of an
isoperimetric region. This fact was unexpected, although this
kind of sets had already appeared as solution in a work by
Pedrosa and Ritor{\'e}~\cite{pr}.

We remark that along this paper, we allow
the Gauss curvature $K$ to be a piecewise continuous function on $M$,
which enlarges the family of surfaces considered.

Although the solutions of the isoperimetric problem in a flat
torus are well-known (they are disks for small values of area, and
bands, see~\cite{H} or ~\cite{hhm1}), the same question in the standard
torus of revolution has remained open for a long time, as
mentioned in~\cite{dp}, an interesting work solving the double
bubble problem in flat tori.
From our classification of stable regions, the solution of
the isoperimetric problem in such surfaces is thus reduced
to numerical comparison between candidates.
In Section~\ref{sec:isoperimetric} we roughly
describe the isoperimetric regions in these surfaces.

We have organized this paper in several sections. In
Section~\ref{sec:preliminares} we establish some notation and
preliminaries, show the constant geodesic curvature curves and
give some stability criteria for them.
Section~\ref{sec:unduloids} is devoted to unduloid type curves:
we mainly prove the existence of closed embedded stable ones in
certain surfaces.
In Section~\ref{sec:stableregions} we study the stability of regions
whose boundary is composed by the curves exposed in
Section~\ref{sec:preliminares}, describing all possible stable
regions in our surfaces (Theorem~\ref{th:main}). Finally in
Section~\ref{sec:isoperimetric}, we make several comments
and future directions of research regarding the isoperimetric problem.

\emph{Acknowledgments.} The author would like to thank Manuel Ritor\'e 
for his continuous support and kind help during the elaboration of
these notes.

\section{Preliminaries}
\label{sec:preliminares}

Let $M$ be a rotationally symmetric torus, that is, a torus
endowed with a one-parameter group of intrinsic isometries. This
kind of surface can be seen as the quotient of a warped product in
the following sense: consider the product $\sph^1 \times I$, where
$I = [-t_0,t_0]$ is a real interval, with the Riemannian metric
\begin{equation}
\label{eq:metrica} ds^2 = dt^2 + f(t)^2 d\theta^2,
\end{equation}
where $\theta\in\sph^1$, $t\in I$, and $f$ is a $C^1$ and
piecewise $C^2$ positive function defined on $I$
(that is, $f''$ is continuous on $I$, except possibly
on a finite number of points). We assume that
$f(t_0)=f(-t_0)$, and $f'(t_0)=f'(-t_0)=0$. Then we can identify
the curves $\sph^1\times\{t_0\}$ and $\sph^1\times\{-t_0\}$ in
order to obtain the torus $M$.

We will also assume that $M$ is symmetric with respect to the
curve $\sph^1\times\{0\}$, that is, $f(t)=f(-t)$ for all $t\in
[-t_0,t_0]$. And that the Gauss curvature $K$ is a decreasing
function of the distance from $\sph^1\times\{0\}$.

The horizontal curves $\sph^1\times \{t\}$ have constant
geodesic curvature and will be called \emph{circles of
revolution} or \emph{parallels}. The vertical
curves $\{\theta\}\times [-t_0,t_0]$ are geodesics of the
metric~\eqref{eq:metrica}. They will be named \emph{vertical
geodesics}.

In this setting, the Gauss curvature $K$ only depends on $t$,
and it is given by
\[
K(t)=-\frac{f''(t)}{f(t)}.
\]
Since $f$ is
only assumed to be a piecewise $C^2$ function, the Gauss curvature
will be a piecewise continuous function in general.

Furthermore, for a circle of revolution $\sph^1 \times \{t\}$, the
length and the geodesic curvature with respect to the normal
vector $-\ptl_t$ are given by
\[
L(t)=2 \pi f(t), \hspace{1cm}      h(t)=\frac{f'(t)}{f(t)}.
\]

An important remark is that the function
\[
(f')^2 - f f''=(2\pi)^{-2}\,L^2(K + h^2)
\]
and $K$ have the same monotone behavior.

Since $K$ is a decreasing function from $\sph^1\times\{0\}$
and $M$ is not a flat torus, we have that $K(0)>0$ and $K(t_0)<0$.
Taking into account that $(f')'=-K f$, we obtain that $f'$
is strictly negative in $(0,t_0)$, so that $f$ is strictly decreasing in
$(0,t_0)$. Hence $\sph^1\times\{0\}$ will be called the
\emph{longest parallel}, and $\sph^1\times\{t_0\}$ will be
the \emph{shortest parallel}. The symmetry of $M$, with
respect to $\sph^1\times\{0\}$, will be referred to as
\emph{horizontal symmetry}.

We now show two examples of the surfaces considered in this work.

\begin{example}
\label{ex:standard}
The standard torus, obtained by rotating a circle
of radius $r$, whose center is at distance $a$ from the axis of
revolution, yields one of these surfaces, with continuous Gauss
curvature. In this case, the interval $I$ can be taken as $[-\pi
r,\pi r]$, and the metric~\eqref{eq:metrica} is given by the
function
\[
f(t)=a+r\,\cos(t/r),
\]
for $t\in I$, where $a>r$.
\end{example}

\begin{example}
\label{ex:tororaro}
Another example is given by the following surface: let $S$ be a
sphere of radius $a$, where two identical disks, centered at each
pole and at height $t^*$ and $-t^{*}$ from the equator,
have been removed. Paste a bounded hyperbolic annulus
of curvature $-b^2$ to each boundary component of $S$
in an appropriate way, in order to
have a $C^1$ and piecewise $C^2$ surface. Finally, by identifying
the two remaining boundary components of the annuli we obtain one
of our surfaces, now with piecewise continuous Gauss curvature.
Here the interval $I$ is $[-d/b,d/b]$ and the
metric~\eqref{eq:metrica} is provided by the function
\[
f(t)=\frac{1}{a}\cos(a\,t),\ t\in [0,t^*],
\]
in the upper spherical piece, and
\[
f(t)=c\,\cosh(d-b\,t),\ t\in [t^*,d/b],
\]
in the upper hyperbolic piece. By the $C^1$-differentiability of
$f$ in $t^*$, it turns that $b$ must be greater than $a\tan(a\,t^*)$
and then
\begin{align*}
c^2&=\frac{1}{a^2}\cos^2(a\,t^*)-\frac{1}{b^2}\sin^2(a\,t^*),
\\
d&=b\,t^* + \cosh^{-1}\bigg(\frac{\cos(a\,t^*)}{a\,c}\bigg).
\end{align*}
\end{example}

\begin{figure}[h]
\centerline{\includegraphics[width=0.45\textwidth]{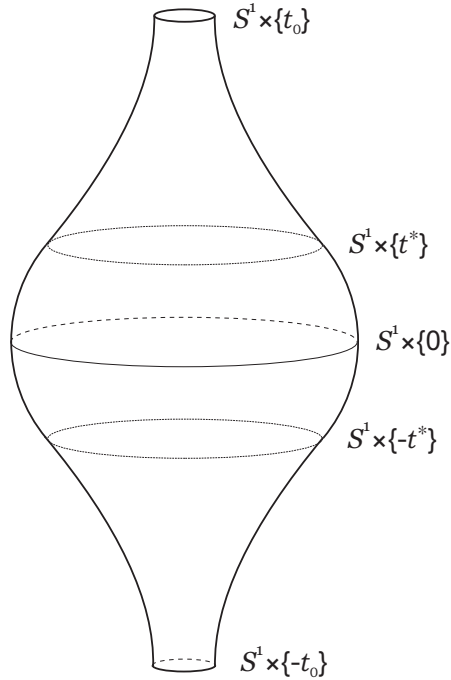}}
\caption{One of the surfaces described in Example~\ref{ex:tororaro},
built from a sphere and two hyperbolic annuli}
\label{fig:standard}
\end{figure}

\begin{remark}
This second example can also be considered with the function
\[
f(t)=\cos(a\,t),\ t\in [0,t^*]
\]
in the upper spherical piece, obtaining other kind of surfaces.
These functions are associated to singular orbifolds with
constant Gauss curvature.
\end{remark}

\subsection{Constant geodesic curvature curves}
\label{subsec:curves}
In this subsection we will describe the closed embedded curves
with constant geodesic curvature appearing in our surfaces.
In the warped product $\sph^1 \times I$, these curves have been
well studied in \cite{ritore} when $f$ is $C^2$, by using a result
by Osserman~\cite[Lemma~7]{O}. However, such a result is also valid
when $f$ is piecewise $C^2$, and so we can extend most of the
consequences to our case.

Let $\ga(s)=(\theta(s),t(s))$ be a curve parametrized
by arc-length $s$ in $S^1\times I$. Denote by $d\ga/ds$ the tangent
vector to $\ga$, and by $\sigma$ the oriented angle
$\angle\,(\ptl_t,d\ga/ds)$. Consider the unit normal vector to $\ga$
given by
\[
N=\frac{\cos\sigma}{f(t)}\,\ptl_{\theta} - \sin\sigma\,\ptl_t,
\]
and let $h$ be the geodesic curvature of $\ga$ with respect to
$N$.

\begin{proposition}$($\cite[Prop.~1.1]{ritore}$)$
\label{prop:ecuaciones}
With the notation above, the curve $\ga$
satisfies the following system of ordinary differential equations
\begin{align}
\label{system} \frac{dt}{ds}&=\, \cos\sigma, \notag
\\
\frac{d\theta}{ds}&=\,\frac{\sin\sigma}{f(t)},
\\
\frac{d\sigma}{ds}&=\,h-\frac{f'(t)}{f(t)}\,\sin\sigma. \notag
\end{align}
Moreover, if h is constant then, for any $c\in I$, the function
\begin{equation}
\label{eq:primeraintegral} f(t)\sin\sigma -
h\,\int_{c}^{t}f(\xi)\,d\xi
\end{equation}
is constant over any solution of \eqref{system}.
\end{proposition}

\begin{remark}
\label{re:periodic}
In view of the system~\eqref{system}, it can be checked that a constant geodesic
curvature curve $\ga$ in $\mathbb{S}^1\times I$ is periodic with respect to
any of the critical points of $t_{|\ga}$.
\end{remark}

The function \eqref{eq:primeraintegral} is called a \emph{first
integral} of \eqref{system}, and allows, as in \cite{ritore},
to classify the closed embedded curves with constant geodesic
curvature in $M$, obtaining the following

\begin{theorem}
\label{th:curves}
Let $M$ be a rotationally symmetric torus with
decreasing Gauss curvature from the longest parallel.
Assume also that $M$ has a horizontal symmetry. Let $C$ be a
connected closed embedded curve in $M$ with constant geodesic
curvature.

Then $C$ is a circle of revolution, a vertical geodesic, a nodoid
type curve or an unduloid type curve.
\end{theorem}

Consider a constant geodesic curvature curve $C$ in $\sph^1\times I$,
with a strict maximum $t(s_0)$ of the $t$-coordinate,
and let $s_1>s_0$ be the next critical point.
If $\sin(\sg(s_1))=1$, the curve will be called unduloid type curve,
and it can be checked that is a periodic graph over $\theta$.
Otherwise, if $\sin(\sg(s_1))=-1$, the curve will be a nodoid type
curve and will present points with vertical tangent vector.

\begin{figure}[ht]
\centerline{\includegraphics[width=0.6\textwidth]{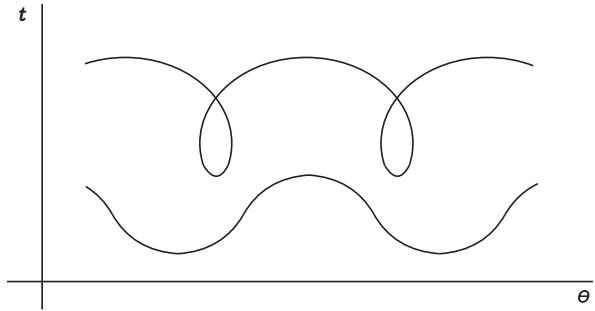}}
\caption{Nodoid and unduloid type curves in $\mathbb{S}^1\times I$}
\label{fig:curvas}
\end{figure}

We shall summarize some properties of circles of revolution in a
warped product $\sph^1\times I$.

\begin{lemma}
\label{le:circulos}
The geodesic curvature of circles of
revolution $h(t)$, computed with respect to the normal vector
$-\ptl_t$, satisfies the following properties:
\begin{itemize}
\item[i)] $h(t)$ is an antisymmetric function on $[-t_0,t_0]$.

\item[ii)] $h(t)$ is increasing in the interval where $(f')^2 - f f''\leq 0$,
and decreasing where
$(f')^2 - f f''\geq 0$.
\end{itemize}
\end{lemma}

\begin{remark}
\label{re:monotonia}
Note that $h(0)=0$, and $h(-t_0)=h(t_0)=0$, since $f'(t_0)=0$.
Moreover, $[(f')^2 - f f''](-t_0)=f(-t_0)^2\,K(-t_0)<0$, and
$[(f')^2 - f f''](0)=f(0)^2\,K(0)>0$.
Then, as $K(t)$ is increasing in $(-t_0,0)$,
Lemma~\ref{le:circulos} gives that $h(t)$ is increasing from $-t_0$
until the point in $(-t_0,0)$ where $(f')^2 - f f''$ vanishes,
and decreasing from that point until $0$.
Consequently, $h(t)$ is positive in $[-t_0,0]$,
and negative in $[0,t_0]$.
\end{remark}

We will now treat when nodoid and unduloid type curves yield closed
embedded curves.
We first define the \emph{period} of an unduloid type curve
as the $\theta$-distance between two
consecutive maxima (or minima) points of the
$t$-coordinate.

\begin{lemma}$($\cite[Prop.~1.3]{ritore}$)$
\label{le:minimax}
Let $C$ be a curve with constant geodesic
curvature in a warped product $\sph^1 \times I$.
\begin{itemize}
\item[i)] If $C$ is a nodoid type curve, it yields a closed
embedded curve if and only if the maximum and the minimum of
$t|_C$ are in the same vertical line.

\item[ii)] If $C$ is an unduloid type curve, it yields a closed
embedded curve if and only if the period of $C$
is equal to $2\pi/k$, with $k\in\nn$.
\end{itemize}
\end{lemma}

\subsection{Stability and the index form}
Consider a curve $C$ with constant geodesic curvature $h$,
not necessarily connected, enclosing a certain area of $M$.
Then, it is well-known that $C$ is a critical point
for the length functional, for area-preserving variations \cite{bapa}.
We shall say that $C$ is \emph{stable} if it is a local minimum of
perimeter for any variation of $C$ with fixed area enclosed.
If $C$ is contained in an open region where $K$ is continuous,
then the second derivative of length is given by
\begin{equation}
\label{eq:indexform} I(u) = -\int_C u\,\bigg\{\frac{d^2 u}{ds^2} +
(K + h^2)\,u\bigg\}ds,
\end{equation}
where $u:C\rightarrow\rr$ is the normal component of the vector
field associated to the variation. In this case, we have
that the stability of $C$ equivalent to
\[
I(u)\geq 0,\ \text{for any function}\ u \ \text{such that}\,
\int_C u\, ds = 0.
\]

A set $\Om\subset M$ is called a \emph{stable region} if $\ptl\Om$
is an embedded stable curve with constant geodesic curvature with
respect to the inner normal. This means that the boundary of a
stable region is a second order local minimum for the perimeter
when keeping constant the area enclosed. It is clear that any
isoperimetric region is stable.

Hereafter, the quadratic form of \eqref{eq:indexform} will be
called the \emph{index form}, and the associated self-adjoint
operator
\begin{equation}
\label{eq:jacobi} J(u)=\frac{d^2 u}{ds^2} + (K + h^2)\,u
\end{equation}
will be called the \emph{Jacobi operator}.

A function $u:C\rightarrow\rr$ satisfying $J(u)=0$ is a
\emph{Jacobi function}.
For instance, the normal component
$u=\escpr{N,\ptl_{\theta}}=f(t)\,\cos\,\sg$
of the Killing vector field $\ptl_{\theta}$
is always a Jacobi function.
Moreover, Jacobi functions also arise from variations of $C$
keeping constant the geodesic curvature along the deformation,
since for a variation with normal component $u$ we have (see~\cite{bgs})
\[
\frac{dh}{dt}\bigg|_{t=0}= u'' + (K+h^2)\,u=J(u).
\]

Given a Jacobi function $u$, a \emph{nodal region} is a connected
component of the complementary in $C$ of the set $\{x\in C:u(x)=0\}$.
By Courant's Nodal Domain Theorem it follows that stable
connected curves have at most two nodal regions (see~\cite[Ch.~I, pag.~19]{chavel}).

The next lemma gives a stability condition for circles of revolution.

\begin{lemma}$($\cite[Lemma~1.6]{ritore}$)$
\label{le:curvatura}
A circle of revolution $\sph^1 \times \{t\}$ is
stable if and only if
\[
(K+h^2)(t)\leq \frac{4\,\pi^2}{L^2(t)},
\quad\text{or equivalently}\quad [(f')^2 - f f''](t)\leq 1.
\]
\end{lemma}

\begin{remark}
Consider $\tilde{t}\in I$ such that
$[(f')^2-f f''](\tilde{t})=1$.
Then the Jacobi operator for the parallel
$\mathbb{S}^1\times\{\tilde{t}\}$ is
\[
J(u)=u''+\frac{1}{f(\tilde{t})^2}\ u.
\]
It is easy to check that $\sin(\theta(s))$ and $\cos(\theta(s))$
are Jacobi functions of $\mathbb{S}^1\times\{\tilde{t}\}$.
Moreover, this parallel is the unique stable one with
two independent Jacobi functions.
\end{remark}

The following lemma treats stable nodoid type curves in $M$. When they are
closed embedded curves, they bound disks in the surface.

\begin{lemma}
\label{le:nodsim} Let $C$ be a closed embedded stable nodoid type
curve in $M$, not contained in a region with constant Gauss
curvature. Then $C$ intersects symmetrically $\sph^1 \times\{0\}$.
\end{lemma}

\begin{proof}
We give an sketch of the proof,
see \cite[Lemmata~2.3 and ~3.4]{ritore} for details.
As $C$ is closed and embedded, an analytical reasoning
\cite[Lemma~7]{O} will give that it cannot be contained
in a region of $M$ with (non-constant) monotone Gauss
curvature. Therefore $C$ will intersect $\mathbb{S}^1\times\{0\}$,
or $\mathbb{S}^1\times\{t_0\}$. By reflecting the curve,
another application of the same reasoning will show
that $C$ is necessarily symmetric with respect to $\mathbb{S}^1\times\{0\}$.
Finally, if $C$ meets $\mathbb{S}^1\times\{t_0\}$,
a suitable function in the index form yield instability.
\end{proof}

\section{Stable unduloid type curves}
\label{sec:unduloids}
In this section we will consider unduloid type curves.
Our aim is to prove that closed embedded stable ones
may occur in certain rotationally symmetric tori.
We remark that these curves did not appear
in any of the surfaces studied in~\cite{ritore}.
Stability will be obtained by using a result from~\cite{hl},
which involves variations of unduloid type curves
by constant geodesic curvature, and requires
the study of the eigenvalue problem associated
to the Jacobi operator.

First, existence of closed embedded unduloid type curves
in some surfaces is guaranteed from the following result.

\begin{lemma}
\label{le:closedunduloids}
Let $M$ be a rotationally symmetric torus
with a horizontal symmetry and decreasing
Gauss curvature $K$ from the longest parallel.
Assume also that there exists a parallel
$\mathbb{S}^1\times\{\tilde{t}\}$ in $M$
satisfying $[(f')^2-f f''](\tilde{t})=1$,
and such that $K$ is smooth and strictly decreasing
in a neighborhood of $\mathbb{S}^1\times\{\tilde{t}\}$.

Then there are closed embedded unduloid type curves
in $M$, close to $\mathbb{S}^1\times\{\tilde{t}\}$.
\end{lemma}

\begin{proof}
Recall that any unduloid type curve is a graph over $\theta$.
Then, writing  $t=t(\theta),\,\sigma=\sigma(\theta)$, the system
of ordinary differential equations obtained from~\eqref{system}
which satisfies such a curve is
\begin{align}
\label{system2}
    \frac{dt}{d\theta}&=\,f(t)\,\cot\,\sigma,
    \\
    \frac{d\sigma}{d\theta}&=\,h\,\frac{f(t)}{\sin\,\sigma} - f'(t).\notag
\end{align}

Let $\ga(\theta,T,h)=(\theta,t(\theta,T,h),\sigma(\theta,T,h))$
denote the solution of \eqref{system2}, with initial conditions
$t(0,T,h)=T$, $\sigma(0,T,h)=\pi/2$ and geodesic curvature $h$.
For $T>\tilde{t}$, this means that $T$ is the maximum value achieved
by $t(\theta,T,h)$.
Call $\tilde{h}=f'(\tilde{t})/f(\tilde{t})$.

Observe that for $T=\tilde{t}$, $h=\tilde{h}$, the solution of
\eqref{system2} is $(\theta,\,\tilde{t},\pi/2)$, the circle of
revolution $\sph^1\times \{\tilde{t}\}$. We are going to see that
for $T$ close enough to $\tilde{t}$, there exist unduloid type
curves with period $2\pi$, so by Lemma~\ref{le:minimax} they
will be closed and embedded.

Note that for $(T,h)$ close enough to $(\tilde{t},\tilde{h})$,
the curves $\ga(\theta,T,h)$ will be unduloid type ones
or circles of revolution by
the first integral~\eqref{eq:primeraintegral},
since $\sin\sg>0$.

Let us define the function $F$ by
\begin{equation}
F(T,h)=\sigma(\pi,T,h).
\end{equation}
Note that if $F(T,h)=\pi/2$, the corresponding unduloid type curve
$\ga(\theta,T,h)$ has period $2\pi$. Clearly $F(\tilde{t},\tilde{h})=\pi/2$.
We want to find functions $h(T)$ such that for $T$ close enough to $\tilde{t}$,
$F(T,h(T))=\pi/2$. Let us compute the partial derivatives of $F$ at $(\tilde{t},\tilde{h})$,
which coincide with the partial derivatives of $\sigma$ at
$(\pi,\tilde{t},\tilde{h})$.

We shall denote by $t_T(\theta)$, $t_h(\theta)$,
(resp. $\sigma_T(\theta)$, $\sigma_h(\theta)),\dots$,
the partial derivatives 
of the function $t(\theta,T,h)$ (resp. $\sigma(\theta,T,h)$)
at $(\theta,\tilde{t},\tilde{h})$. Then,
Taylor developments at $(\theta,\,\tilde{t},\,\tilde{h})$ give
\begin{align*}
t(\theta,T,h)&=\tilde{t}+(T-\tilde{t})\,t_T(\theta) + (h-\tilde{h})\,t_h(\theta) +
(T-\tilde{t})(h-\tilde{h})\,t_{Th}(\theta) +
\\
&\frac{1}{2}(T-\tilde{t})^2\,t_{TT}(\theta) +
\frac{1}{2}(h-\tilde{h})^2\,t_{hh}(\theta) + \frac{1}{6} (T-\tilde{t})^3\,t_{TTT}(\theta)
+ \dots,
\\
\sigma(\theta,T,h)&=\pi/2 + (T-\tilde{t})\,\sigma_T(\theta) + (h-\tilde{h})\,\sigma_h(\theta) +
(T-\tilde{t})(h-\tilde{h})\,\sigma_{Th}(\theta) +
\\
&\frac{1}{2}(T-\tilde{t})^2\,\sigma_{TT}(\theta) +
\frac{1}{2}(h-\tilde{h})^2\,\sigma_{hh}(\theta) + \frac{1}{6} (T-\tilde{t})^3\,\sigma_{TTT}(\theta)+
\dots.
\end{align*}
Moreover, from the definitions of $t(\theta,T,h)$ and $\sigma(\theta,T,h)$,
we have that $t_T(0)=1$, $t_h(0)=t_{Th}(0)=t_{TT}(0)=t_{hh}(0)=t_{TTT}(0)=0$,
$\sigma_T(0)=\sigma_h(0)=\sigma_{Th}(0)=\sigma_{TT}(0)=\sigma_{hh}(0)=\sigma_{TTT}(0)=0$.

From \eqref{system2}, considering the second equation as
\[
\frac{d\sigma}{d\theta}=\,(h-\tilde{h})\,\frac{f(t)}{\sin\,\sigma} + \tilde{h}\,\frac{f(t)}{\sin\,\sigma} - f'(t),
\]
and Taylor developments of the involved functions, we get
\begin{align*}
\frac{dt_T}{d\theta}&=-f(\tilde{t})\,\sigma_T,
\\
\frac{d\sigma_T}{d\theta}&=\frac{t_T}{f(\tilde{t})},
\end{align*}
which gives $t_T(\theta)=\cos(\theta)$, $\sigma_T(\theta)=\sin(\theta)/f(\tilde{t})$.

In the same way,
\begin{align*}
\frac{dt_h}{d\theta}&=-f(\tilde{t})\,\sigma_h,
\\
\frac{d\sigma_h}{d\theta}&=f(\tilde{t}) + \frac{t_h}{f(\tilde{t})},
\end{align*}
obtaining $t_h(\theta)=f(\tilde{t})^2\,(\cos(\theta)-1)$, $\sigma_h(\theta)=f(\tilde{t})\,\sin(\theta)$.

Therefore, the gradient of $F$ at $(\tilde{t},\tilde{h})$ is equal to
\[
\nabla F(\tilde{t},\tilde{h})=(\sigma_T(\pi),\sigma_h(\pi))=(0,0).
\]

We now compute the hessian of $F$ at $(\tilde{t},\tilde{h})$.
From the equations above we obtain the following systems.
\begin{align}
\frac{dt_{TT}}{d\theta}&=-f(\tilde{t})\,\bigg(\sigma_{TT} - 2\,\tilde{h}\,t_T\,\sigma_T\bigg),
\\
\frac{d\sigma_{TT}}{d\theta}&=f(\tilde{t})\,\bigg(\frac{1}{f(\tilde{t})^2}\,t_{TT} + K'(\tilde{t})\,t_T^2
+ \tilde{h}\,\sigma_T^2\bigg)\notag,
\end{align}

\begin{align}
\frac{dt_{hh}}{d\theta}&=-f(\tilde{t})\,\bigg(\sigma_{hh} + 2\,\tilde{h}\,t_h\,\sigma_h\bigg),
\\
\frac{d\sigma_{hh}}{d\theta}&=f(\tilde{t})\,\bigg(2\,\tilde{h}\,t_h + \frac{1}{f(\tilde{t})^2}\,t_{hh}
+ K'(\tilde{t})\,t_h^2 + \tilde{h}\,\sigma_h^2\bigg)\notag,
\end{align}

\begin{align}
\frac{dt_{Th}}{d\theta}&=-f(\tilde{t})\,\bigg(\sigma_{Th} + \tilde{h}\,(t_T\,\sigma_h + t_h\,\sigma_T)\bigg),
\\
\frac{d\sigma_{Th}}{d\theta}&=f(\tilde{t})\,\bigg(\tilde{h}\,t_T + \frac{1}{f(\tilde{t})^2}\,t_{Th}
+ K'(\tilde{t})\,t_h\,t_T + \tilde{h}\,\sigma_T\,\sigma_h\bigg)\notag.
\end{align}

Solving these systems, we get that $\sigma_{TT}(\pi)=0$, $\sigma_{Th}(\pi)=\rho/2$ and
$\sigma_{hh}(\pi)=\rho\,f(\tilde{t})^2$, where $\rho=\pi\,f(\tilde{t})^3\,K'(\tilde{t})$,
and so the hessian of $F$ at $(\tilde{t},\tilde{h})$
is given by
\begin{equation}
\label{eq:hessian}
\nabla^2 F(\tilde{t},\tilde{h})= \left( \begin{array}{cc}
0 & \rho/2 \\
\rho/2 & \rho\,f(\tilde{t})^2
\end{array} \right).
\end{equation}
Therefore $\nabla^2 F(\tilde{t},\tilde{h})$ is non-degenerate,
since $K$ is a strictly decreasing function.

Then $(\tilde{t},\tilde{h})$ is a non-degenerate
critical point of $F$. By applying Morse's Lemma \cite[Lemma~2.2]{mi}
we obtain two planar curves $\alpha_c,\alpha_o$ with
$\alpha_c(0)=\alpha_o(0)=(\tilde{t},\tilde{h})$ and such that,
for any point $(T,h)$ close to $(\tilde{t},\tilde{h})$ satisfying
$F(T,h)=\pi/2$, it follows that $(T,h)$ lies on the trace
of $\alpha_c\cup\alpha_o$.
Let us study the tangent vectors of both curves at the origin.
Denoting by $\alpha$ one of those curves, we know that $F\circ\alpha=\pi/2$.
Differentiating such an equality we obtain $\escpr{\nabla F_\alpha,\alpha'}=0$.
Differentiating once again and evaluating at the origin, we have
\[
\alpha'\,\escpr{\nabla F_{(\tilde{t},\tilde{h})},\alpha'}=
\nabla^2F_{(\tilde{t},\tilde{h})}\,(\alpha',\alpha')=0.
\]
From \eqref{eq:hessian}, we get $\alpha_c'(0)=(1,-1/f(\tilde{t})^2)$
and $\alpha_o'(0)=(1,0)$. Hence both curves can be written in
terms of $T$ near $\tilde{t}$, because the tangent vectors are not vertical.
So there exist two functions $h_c(T),h_o(T)$ defined in a neighborhood of
$\tilde{t}$ such that $\{(T,h_c(T))\}$ is the trace of the curve $\alpha_c$,
and $\{(T,h_o(T))\}$ is the trace of $\alpha_o$. Initial conditions $(T,h_c(T))$
in \eqref{system2} yield circles of revolution,
and $(T,h_o(T))$ give unduloid type curves as solutions (note that
$h_c'(\tilde{t})=-1/f(\tilde{t})^2$ and $h_o'(\tilde{t})=0$).

Since
\begin{equation}
\label{eq:pi/2}
F(T,h_o(T))=(F\circ\alpha_o)(s)=\pi/2,
\end{equation}
we conclude that the unduloid type curves
$\ga(\theta,T,h_o(T))$ have period $2\pi$,
and so they are closed and embedded.
\end{proof}

We now focus on the stability of closed embedded unduloid type curves.
A first remark is that stable ones will present
a \emph{unique} maximum point (and therefore a unique minimum point)
for the $t$-coordinate. Otherwise, the normal component of the
rotations vector field will have more than two nodal regions,
yielding instability by applying Courant's Nodal Domain Theorem.
In view of Lemma~\ref{le:minimax}, this fact is equivalent to
that the period of closed embedded stable unduloid type curves
equals $2\pi$.

\begin{figure}[ht]
\centerline{\includegraphics[width=0.45\textwidth]{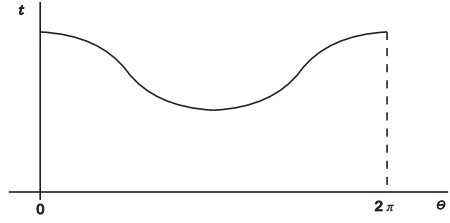}}
\caption{A closed embedded stable unduloid type curve in $\mathbb{S}^1\times I$}
\label{fig:stableunduloid}
\end{figure}

The following lemma states another necessary condition for stability.

\begin{lemma}
\label{le:unduloids} Let $C$ be a closed embedded stable unduloid
type curve in $M$, not contained in the region $\{(f')^2 - f f'' = 1\}$.
Then there are points of $C$ satisfying $(f')^2 - f f'' > 1$, and other
ones satifying $\,(f')^2 - f f'' < 1$.
\end{lemma}

\begin{proof}
Suppose first that $C$ lies in a region where $(f')^2 - f f'' \leq 1$.
Then, by \cite[Lemma~2.3]{ritore}, its period is greater than
$2\pi$, and so it is not a closed embedded curve by
Lemma~\ref{le:minimax}.
Assume now that $C$ lies in a region where $(f')^2 - f f'' \geq 1$.
Then its period is less than $2\pi$ (see~\cite[Lemma~2.13]{ritore}
for details), and so $C$ will present more than one maximum point
for the $t$-coordinate, contradicting stability.
\end{proof}

In view of the Jacobi operator~\eqref{eq:jacobi},
we can consider the eigenvalue problem associated to a
closed embedded unduloid type curve $C$ in $M$:
\begin{equation}
\label{eq:eigenproblem}
J(u) + \lambda\,u=0,
\end{equation}
with $u:C\to\rr$ a $C^2$ function. 

\noindent We will call eigenvalues to the real numbers $\lambda$ for
which there exist functions satisfying \eqref{eq:eigenproblem}.
Such functions will be named eigenfunctions associated to $\lambda$.
We recall some well-known general facts about this problem.

\begin{lemma}$($\cite[Chapter~8, Theorem~2.1]{cl},\cite{chavel}$)$
\label{le:eigen}
Given a curve $C$, the eigenvalues associated to the Jacobi operator
form an increasing sequence $\{\lambda_i\}_{i\geq 1}$.
Furthermore, the space of eigenfunctions $V_{\lambda_i}$
associated to $\lambda_i$, when considering an eigenvalue problem
with boundary conditions, is one-dimensional
and each $\phi_i\in V_{\lambda_i}$ has exactly $i-1$ zeros.
\end{lemma}

We will now show some results regarding
the eigenvalue problem~\eqref{eq:eigenproblem},
necessary for our purposes.

\begin{lemma}
\label{le:valorpropio}
Let $C$ be a closed embedded unduloid type curve in $M$.
Then the first eigenvalue for the Jacobi operator
associated to $C$ is negative,
and the second eigenvalue is non-positive.
\end{lemma}

\begin{proof}
Consider the normal component $u$ of the rotations vector field on $M$,
which is a Jacobi function.
Therefore the restriction of $u$ to $C$ is
an eigenfunction for the associated eigenvalue problem~\eqref{eq:eigenproblem}
with zero eigenvalue.

It is easy to check that $u$ gives at least two nodal regions in $C$,
since the period of $C$ must be $2\pi/k$, for some $k\in\nn$,
by Lemma~\ref{le:minimax}.
Therefore, by applying Courant's Nodal Domain Theorem, the statement follows.
\end{proof}

\begin{remark}
\label{re:cero}
Closed embedded stable curves have
non-negative second eigenvalue $\lambda_2$
for the above eigenvalue problem~\eqref{eq:eigenproblem}.
From Lemma~\ref{le:valorpropio} it follows that
closed embedded stable unduloid type curves in $M$
verify $\lambda_2=0$.
\end{remark}

\begin{remark}
\label{re:eigennodoid}
Given $C$ a nodoid type curve in $M$, the same reasoning as above
shows that Lemma~\ref{le:valorpropio} also holds
when $C$ is closed and embedded.
\end{remark}

Let $\lambda$ be an eigenvalue of~\eqref{eq:eigenproblem} for a
given closed embedded unduloid type curve $C$,
with associated eigenfunction $u$.
Call \emph{fundamental piece of $C$} to any region of $C$
delimited by a maximum point of the $t$-coordinate, and
the consecutive minimum point.
By considering the vertical reflection with respect to a
maximum point of the $t$-coordinate of $C$, we
can express $u=u_s + u_a$, where $u_s$ is a symmetric eigenfunction
satisfying the Neumann boundary condition in any fundamental piece $C'$ of $C$,
and $u_a$ is an antisymmetric eigenfunction satisfying the Dirichlet
boundary condition in $C'$, both of them with eigenvalue $\lambda$.
We will also name $\lambda_i^N(C')$, $\lambda_i^D(C')$
to the eigenvalues associated to $C'$ for the Neumann and Dirichlet eigenvalue
problem, respectively.

Recall that $\lambda_1<0$ by Lemma~\ref{le:valorpropio}.
Therefore, the above reasoning leads to a negative Dirichlet
eigenvalue in $C'$, which is not possible since $\lambda_1^D(C')=0$.
Then $u_a=0$, and so $u=u_s$. Consequently, $\lambda_1$ coincides with
$\lambda_1^N(C')$.

Now consider the second eigenvalue $\lambda_2$ in $C$.
If it is negative, then the curve is unstable, and moreover,
an analogous treatment will give
$\lambda_2=\lambda_2^N(C')$~\cite{bb}.
On the other hand, if
$\lambda_2\geq 0$, then necessarily $\lambda_2=0$
by Lemma~\ref{le:valorpropio}.

We will now study the period of unduloid type curves,
not necessarily closed and embedded.
For $\ga(\theta,T_0,h)$, such a period is defined as
the $\theta$-distance between two consecutive
maximum points (or minimum points) for the $t$-coordinate.
If we move slightly the maximum point $T$,
keeping constant the geodesic curvature,
the period of $\ga(\theta,T,h)$ only
depends on $T$, and then the derivative of the period
with respect to $T$, at $T=T_0$,
can be considered.

\begin{lemma}$($\cite[Corollary~2.8]{pr}$)$
\label{le:Neumann}
Consider an unduloid type curve $C$ in M,
close enough to $\mathbb{S}^1\times\{\tilde{t}\}$.
Let $C'$ be a fundamental piece of $C$. Then the second
eigenvalue $\lambda_2^{N}(C')$ of the Neumann problem for
the Jacobi operator in $C'$ is positive
if and only if the derivative of the period with respect to
the maximum point is positive.
\end{lemma}

\begin{proof}
Call $C=\ga(\theta,T_0,h)$ and consider the variation by
unduloid type curves given by $\ga(\theta,T,h)$,
with $T$ close to $T_0$ and $h$ fixed.
It can be checked that the associated variational function
$u$ is a Jacobi function satisfying $u(0)=1$, $u'(0)=0$ and $u''(0)<0$,
and so $u'(\varepsilon)<0$ for positive $\varepsilon$ close to $0$.

Assume first that the derivative is positive,
and let us find a convenient expression of it.
Let $\theta_2>0$ be the first instant where
the $t$-coordinate of $C$ achieves a minimum.
Therefore $\sg'(\theta_2)<0$,
and $[0,\theta_2]$ yields a fundamental piece
of $C$, say $C'$.
By applying Implicit Function Theorem to $\sg(\theta,T)$
at $(\theta_2,T_0)$,
we obtain a function $\theta(T)$ such that
$\sg(\theta(T),T)=\pi/2$, and so $\theta(T)$
gives (half) the period of $\ga(\theta,T,h)$.
Straightforward calculations for computing $\theta'(T)$
show that the desired derivative equals
\[
\frac{2\,u'(\theta_2)}{f(t(\theta_2,T_0))\ \sg'(\theta_2)}.
\]
Since we are assuming that such a value is positive,
then we have $u'(\theta_2)<0$.

We will focus now on the interval $(0,\theta_2)$.
Observe that $u$ will have at most one zero on such interval.
Otherwise, by considering two consecutive zeroes $\theta_3$,
$\theta_4$, it turns that the first eigenvalue of
the Dirichlet problem in $(\theta_3,\theta_4)$ is zero
(recall that $u$ is a Jacobi function),
and then $\lambda_1^D(C')$ will be strictly negative,
by the monotone property of eigenvalues \cite{chavel},
which is contradictory.

On the other hand, $u'$ does not vanish in $(0,\theta_2)$.
Since $u'(\varepsilon)<0$ for $\varepsilon>0$,
and $u'(\theta_2)<0$, if $u'$ has a zero
then $u''$ will have two zeroes.
As $u''+(K+h^2)\,u=0$, it follows that $u$
will vanish twice in $(0,\theta_2)$,
and we would proceed as above.

Let $\theta_3$ be the first zero of $u'$
greater than $\theta_2$.
As $u'$ is strictly negative in $(0,\theta_3)$
and vanishes at the extremes,
then $u''$ changes its sign in $(0,\theta_3)$.
It follows that $u$ is strictly decreasing and
vanishes only once in $(0,\theta_3)$.
By Lemma~\ref{le:eigen},
$u$ is the second eigenfunction of the Neumann problem
in $C|_{(0,\theta_3)}$,
and the associated second Neumann eigenvalue is zero.
By the monotone property of eigenvalues we conclude that
$\lambda_2^{N}(C')>0$.

Assume now that $\lambda_2^{N}(C')>0$.
If the derivative of the period is negative,
by applying directly \cite[Cor.~2.8]{pr} we
obtain that $\lambda_2^{N}(C')<0$, a contradiction.
If the derivative of the period is zero,
same reasoning as above leads to
$\lambda_2^{N}(C')=0$, again a contradiction.
So the derivative must be positive
and the statement follows.
\end{proof}

The next lemma determines the sign of the derivative
of the period of closed embedded unduloid type curves
(close enough to $\mathbb{S}^1\times\{\tilde{t}\}$)
by means of the following intrinsic condition
on the surface.

\begin{lemma}
\label{le:derivadaperiodo}
Consider a surface $\mathbb{S}^1\times I$ under the conditions of
Lemma~\ref{le:closedunduloids}.

Then, the derivative of the period of
closed embedded unduloid type curves
with respect to $T$ is strictly positive,
for $T>\tilde{t}$ close enough,
if and only if
\begin{equation}
\label{eq:condicion}
3\,K(\tilde{t})\,(1-f(\tilde{t}))\,+
3\,f(\tilde{t})^2\,(\tilde{h}\,K'(\tilde{t})- K''(\tilde{t})) +
5\,f(\tilde{t})^4\,K'(\tilde{t})^2 >0.
\end{equation}
\end{lemma}

\begin{proof}
Fix $T_0>\tilde{t}$ close enough, and $h_o(T_0)$ as geodesic curvature
in order to compute the desired derivative at $T=T_0$.
It is clear that such a derivative will be strictly positive
if and only if  the derivative of the period at $T=\tilde{t}$
is strictly positive, for constant geodesic curvature
$\tilde{h}=h(\tilde{t})$
(note that the derivative will preserve the same sign
in an approppiate neighborhood of $(\tilde{t},\tilde{h})$).

Denote by $t(\theta,T)$, $\sigma(\theta,T)$ to
the solutions of \eqref{system2} with initial
conditions $t(0,T)=T$, $\sigma(0,T)=\pi/2$, and
geodesic curvature $\tilde{h}$.

We search a function $\theta:\rr\to\rr$ verifying
$\sigma(\theta(T),T)=\pi/2$. If such a function exists,
$\theta(T)$ will give the (half) period of $\ga(\theta,T,\tilde{h})$.

Consider the auxiliar function $F$ given by
\begin{equation}
\label{eq:defF}
\sigma(\theta,\,T)-\pi/2=(T-\tilde{t})\,F(\theta,T),
\end{equation}
and extended continuously in $(\theta,\tilde{t})$.
That is,
\[
F(\theta,\tilde{t})=\sg_T(\theta,\tilde{t})=\frac{\sin(\theta)}{f(\tilde{t})}.
\]
Observe that $F(\pi,\tilde{t})=0$. Moreover,
for $T\neq\tilde{t}$, we have that $F(\theta,T)=0$
if and only if $\sg(\theta,T)=\pi/2$.

From Taylor development of $\sigma(\theta,\,T)$
of Lemma~\ref{le:closedunduloids} we have
\begin{equation}
\label{eq:defF2}
F(\theta,T)=\sigma_T(\theta) + \frac{1}{2}\,\sigma_{TT}(\theta)\,(T-\tilde{t}) +\frac{1}{6}\,\sg_{TTT}(\theta)\,(T-\tilde{t})^2 +  O((T-\tilde{t})^3).
\end{equation}
By differentiating \eqref{eq:defF2},
\[
\frac{\ptl F}{\ptl\theta}(\pi,\tilde{t})=\sigma_T'(\pi)=\frac{-1}{f(\tilde{t})}\neq 0.
\]
By applying the Implicit Function Theorem, there exists $\theta(T)$,
for $T$ close to $\tilde{t}$, such that $\theta(\tilde{t})=\pi$,
and
\begin{equation}
\label{eq:F}
F(\theta(T),T)=0,
\end{equation}
equivalently $\sigma(\theta(T),T)=\pi/2$.
So this function $\theta(T)$ gives the period
of the unduloid type curve $\ga(\theta,T,\tilde{h})$.

We now compute the derivative of $\theta(T)$.
Since
\[
\theta(T)=\theta(\tilde{t}) + \theta'(\tilde{t})\,(T-\tilde{t}) +
\frac{1}{2}\,\theta''(\tilde{t})\,(T-\tilde{t})^2 + O((T-\tilde{t})^3),
\]
we have
\begin{equation}
\label{eq:derivadaperiodo}
\theta'(T)=\theta'(\tilde{t}) + \theta''(\tilde{t})\,(T-\tilde{t}) +
O((T-\tilde{t})^2).
\end{equation}

From \eqref{eq:F} it follows that
\[
\theta'(T)=-\frac{\ptl F/\ptl T}{\ptl F/\ptl\theta}(\theta(T),T).
\]
Evaluating in $T=\tilde{t}$, and taking into
account \eqref{eq:defF2}, we have
\[
\theta'(\tilde{t})=\frac{-\sigma_{TT}(\pi)}{2\,\sigma_T(\pi)}=
\frac{1}{2}f(\tilde{t})\,\sigma_{TT}(\pi).
\]
From \eqref{eq:hessian}, $\sigma_{TT}(\pi)=0$, and hence $\theta'(\tilde{t})=0$.

Differentiating once again \eqref{eq:defF2} and evaluating at $T=\tilde{t}$,
since $\theta'(\tilde{t})=0$ we get
\[
\theta''(\tilde{t})=\frac{-\ptl^2 F/\ptl T^2}{\ptl F/\ptl\theta}(\pi,\tilde{t})=
\frac{-\sigma_{TTT}(\pi)}{3\,\sg_T'(\pi)}=\frac{1}{3}f(\tilde{t})\,\sg_{TTT}(\pi).
\]

In view of \eqref{eq:derivadaperiodo}, the sign of $\theta'(T)$,
for $T>\tilde{t}$, depends on $\sg_{TTT}(\pi)$.
In order to compute this value, we first have to solve
the system of differential equations for $t_{TTT}$ and $\sigma_{TTT}$,
which can be obtained as in Lemma~\ref{le:closedunduloids}.
Straightforwards calculations give
\begin{align}
\label{sigma3}
\sigma_{TTT}(\pi)=\frac{\pi}{8\,f(\tilde{t})}\,
\bigg\{&3\,K(\tilde{t})\,(1-f(\tilde{t}))\,+
\\
&3\,f(\tilde{t})^2\,(\tilde{h}\,K'(\tilde{t})- K''(\tilde{t})) +
5\,f(\tilde{t})^4\,K'(\tilde{t})^2\bigg\}, \notag
\end{align}
which finishes the proof.
\end{proof}

We finally state the key result in order to prove
the existence of closed embedded unduloid type curves
which are stable.
We will apply the results appearing in~\cite{hl}
(see also \cite{koiso}).
Consider a torus of revolution under the conditions
of Lemma~\ref{le:closedunduloids},
and let $\ga_{T_0}=\ga(\theta,T_0,h_o(T_0))$ be
a closed embedded unduloid type curve, with $T_0$
close enough to $\tilde{t}$, $T_0>\tilde{t}$.
We know that $\ga(\theta,T,h_o(T))$ is a family
of closed embedded unduloid type curves, for $T$ close to $T_0$.
Then we have
\begin{lemma}$($\cite[Lemma~2]{hl},\,\cite[Theorem~1.3]{koiso}$)$
\label{le:hl}
Assume $\lambda_1<0\leq\lambda_2$ for the Jacobi operator in $\ga_{T_0}$,
with the notation above. Then $\ga_{T_0}$ is stable if and only if
\[
\frac{dh_o}{dT}\,\frac{da}{dT}<0,
\]
at $T=T_0$, where $\frac{da}{dT}$ is the change of
area induced by the variation $\ga(\theta,T,h_o(T))$.
\end{lemma}

\subsection{Standard tori of revolution}
We now focus on standard tori of revolution
of Example~\ref{ex:standard}, given by
$f(t)=a + r\cos(t/r)$ and $\tilde{t}=\pi r/2$.
We will see that there exist closed embedded
unduloid type curves which are stable
in some of these surfaces.

Firstly observe that the assumptions
of Lemma~\ref{le:closedunduloids}
are verified on any standard torus,
and so there will exist a family of
closed embedded unduloid type curves
parametrized by the maximum point of the $t$-coordinate.

We also point out that the condition shown
in Lemma~\ref{le:derivadaperiodo} is satisfied,
since the left term in \eqref{eq:condicion} is equal to
\[
\frac{5\,a^2 + 9\,r^2}{r^4},
\]
which is strictly positive.
Therefore the derivative of the period of
unduloid type curves with respect to the maximum point
will be positive, for $T>\tilde{t}$.
This allows to compute the second eigenvalue
for the Jacobi operator associated to these curves.

\begin{lemma}
\label{le:lambda2}
Let $C$ be a closed embedded unduloid type curve
in a standard torus of revolution, close enough to
$\mathbb{S}^1\times\{\tilde{t}\}$.
Then the second eigenvalue $\lambda_2$ for the
Jacobi operator in $C$ is equal to zero.
\end{lemma}

\begin{proof}
From Lemma~\ref{le:valorpropio} we know that $\lambda_2\leq 0$.
Assume that $\lambda_2<0$.
Hence, reasonings of this Section give that,
for any fundamental piece $C'$ of $C$,
$\lambda_2^{N}(C')=\lambda_2<0$. This fact
yields a contradiction with Lemma~\ref{le:Neumann},
and so $\lambda_2=0$.
\end{proof}

\begin{lemma}
\label{le:stableunduloids}
Let $M$ be a standard torus of revolution
of Example~\ref{ex:standard} with $r<a<3\,r$.
Then, there exist closed embedded
unduloid type curves which are stable.
\end{lemma}

\begin{proof}
Consider the family $\ga(\theta,T,h_o(T))$ of
closed embedded unduloid type curves given in
Lemma~\ref{le:closedunduloids},
with $T>\tilde{t}$ close enough.
Fix $C$ one of these curves.
By Lemmae~\ref{le:valorpropio} and~\ref{le:lambda2},
we can apply Lemma~\ref{le:hl} to study the stability of $C$.

By differentiating \eqref{eq:pi/2}, it follows that
\[
F_T + F_h\,h_o'=0
\]
and
\[
F_{TT} + 2\,F_{Th}\,h_o' + F_{hh}\,(h_o')^2 + F_h\,h_o''=0.
\]

Differentiating once again and evaluating
at $T=\tilde{t}$, taking into account that
$h_o'(\tilde{t})=0$ and
$F_h(\tilde{t},\tilde{h})=\sigma_h(\pi)=0$, we obtain
\[
F_{TTT}(\tilde{t},\tilde{h}) + 3\,F_{Th}(\tilde{t},\tilde{h})\,h_o''(\tilde{t})=0,
\]
and so
\begin{equation}
\label{eq:ho2}
h_o''(\tilde{t})=\frac{-F_{TTT}(\tilde{t},\tilde{h})}{3\,F_{Th}(\tilde{t},\tilde{h})}=
\frac{-\sigma_{TTT}(\pi)}{3\,\sigma_{Th}(\pi)}=\frac{-\sigma_{TTT}(\pi)}{3\,\rho/2},
\end{equation}
with $\rho$ as defined in~\eqref{eq:hessian}. Since $M$ is a standard torus of
revolution, by~\eqref{sigma3} we conclude that
\[
h_o''(\tilde{t})=\frac{9\,r^2+5\,a^2}{12\,a^3\,r^2}>0,
\]
so $h_o'(T)$ is strictly increasing in $T=\tilde{t}$, and
then positive for $T>\tilde{t}$.

On the other hand, since the associated vector field
induced by the variation $\ga(\theta,T,h_o(T))$ is $(t_T + h_o'\,t_h)\,\ptl_t$,
we have that the derivative of the area along this deformation
is equal to
\[
\int_0^{2\pi} f(T) \bigg(t_T(\theta) + h_o'(T)\,t_h(\theta)\bigg)d\theta,
\]
which vanishes when evaluated at $T=\tilde{t}$.

Moreover, the second derivative of the area
at $T=\tilde{t}$ is given by
\[
\int_0^{2\pi}(f'(\tilde{t})\,t_T(\theta)^2 + f(\tilde{t})\,(t_{TT}(\theta) + h_o''(\tilde{t})\,t_h(\theta)))\,d\theta,
\]
which equals
\[
\frac{(a^2 - 9\,r^2)\,\pi}{6\,r^2}<0,
\]
since $a<3\,r$. Hence, the derivative of the area is
strictly decreasing in $T=\tilde{t}$,
and so, strictly negative for $T>\tilde{t}$. By applying
Lemma~\ref{le:hl}, we conclude that $C$ is stable.
\end{proof}

\section{Stable regions}
\label{sec:stableregions}

In this section we will describe the different stable regions that can
appear in our surfaces. As we know which are the stable constant geodesic
curvature curves by the previous Section, we will
check which combinations of them bound stable regions.

\begin{lemma}$($\cite[Lemma~1.7]{ritore}$)$
\label{le:stabannu} Consider an annulus $\sph^1 \times
[t_1,t_2]$ in $M$. Its boundary is stable if and only if each
component $\sph^1 \times \{t_i\}$ is stable $(i=1,2)$ and
\begin{equation}
\label{eq:nosim}
\frac{K + h^2}{L}(t_1) + \frac{K + h^2}{L}(t_2)\leq 0.
\end{equation}
\vspace{0.2cm} In the case of a symmetric annulus $\sph^1 \times
[-t,\, t]$, above conditions are equivalent to
\[
(K + h^2)(t)\leq 0.
\]
\end{lemma}

The above result completely characterizes the stable annuli bounded by
parallels. It yields that $\sph^1 \times [t,-t]$ is stable for
$t\in (-t_0,t_c)$, with $t_c$ the supremum of the points
in $[-t_0,0]$ where $K+h^2$ (or equivalently $(f')^2 - f f''$)
is nonpositive. Moreover, $\sph^1\times [t_c,-t_c]$ is stable
if and only if $(K+h^2)(t_c)=0$ (recall the possible discontinuities of $K$).
As $h'=-(K+h^2)$, it is clear that $h(t)$ is increasing in $(-t_0,t_c)$
and decreasing in $(t_c,0)$.
Let $\tilde{t}\in [t_c,0]$ be the infimum ot the points
where $L^2 (K+h^2)\geq 4\,\pi^2$ (equivalently $(f')^2 - f f''\geq 1$).
Then, taking into account Lemma~\ref{le:curvatura} and Remark~\ref{re:monotonia},
it follows that for any $t'\in (t_c,\tilde{t})$, $\mathbb{S}^1\times\{t'\}$
is stable and there exists a unique $t''\in (-t_c,t_0)$
satisfying that $\sph^1\times [t',t'']$ is a
nonsymmetric annulus whose boundary has constant geodesic
curvature.

In the light of Lemma~\ref{le:unduloids}, stable domains bounded by an
unduloid type curve and a circle of revolution will appear
if all nonsymmetric annuli $\mathbb{S}^1\times[t',t'']$,
for $t'\in(t_c,\tilde{t})$, are stable,
and the first of such domains will arise
from the annulus corresponding to $t'=\tilde{t}$.

A vertical annulus is a set bounded by two vertical geodesics.
The following result proves that the union of vertical annuli is
always stable.
\begin{lemma}
\label{le:av} Any union of finite disjoint vertical annuli is a stable
region.
\end{lemma}

\begin{proof}
Given a vertical geodesic $C$, we shall first prove that its
length is less than or equal to the length of any other closed
embedded curve $D$, with the same type of homothopy
and near enough to $C$.

Let $N$ denote the unit normal vector field to the set of all
vertical geodesics, which can be properly extended to $M$.
In fact, $N=\frac{1}{f}\,\ptl\theta$. Then, $div\,(N)=0$. We can
assume without loss of generality that $D$ intersects $C$,
consider two consecutive intersection points $p_1,\,p_2$, and call
$\Sigma$ to one of the domains bounded by the pieces of $C$ and $D$ between
$p_1$ and $p_2$. By applying the Divergence Theorem to $N$ in
$\Sigma$, we easily conclude that the length of $C$ is less than
or equal to the length of $D$.

Now it is clear that for a vertical annulus $\Om$, any variation
of $\ptl\Om$ preserving the area enclosed will give more perimeter at
each instant, so $\ptl\Om$ will be a local minimum for the length
while keeping constant the area. Therefore $\Om$ is stable. Notice
that the same argument holds for an arbitrary union of disjoint vertical
annuli.
\end{proof}

Now we can state our Main Theorem, describing all possible stable regions
in our surfaces.

\begin{theorem}
\label{th:main}

Let $M$ be a rotationally symmetric torus with a horizontal
symmetry and with possibly discontinuous, decreasing Gauss
curvature from the longest parallel. Then the stable
regions in $M$ may be:
\begin{itemize}
\item[i)] disks bounded by constant geodesic curvature curves,
which are symmetric with respect to the
longest parallel, or contained in a region with constant
Gauss curvature, and their complements,

\item[ii)] annuli symmetric with respect to the shortest parallel,
bounded by circles of revolution contained in the region $K + h^2
\leq 0$, and their complements,

\item[iii)] nonsymmetric annuli bounded by circles of revolution
contained in the region $K + h^2 \leq\frac{4\,\pi^2}{L^2}$ and
verifying condition~\eqref{eq:nosim}, and their complements,

\item[iv)] unions of vertical annuli bounded by vertical geodesics,

\item[v)] annuli bounded by an unduloid type
curve satisfying Lemma~\ref{le:unduloids}, and a circle of
revolution contained in $K + h^2 < 0$, and their complements,

\item[vi)] unions of a disk and a symmetric annulus with the
same geodesic curvature, in
the above conditions, and their complements.

\end{itemize}
\end{theorem}

\begin{proof}
Let $\Om$ be a stable region in $M$. Then $\ptl\Om$ is an
embedded stable curve (not necessarily connected) with constant
geodesic curvature with respect to the inner normal. By
Theorem~\ref{th:curves}, $\ptl\Om$ will consist of a union of
circles of revolution, vertical geodesics, nodoid type curves and
unduloid type curves.

First observe that no more than two circles of revolution can
appear in $\ptl\Om$, because the geodesic curvatures will not
coincide in view of Remark~\ref{re:monotonia}.

Let $C_1$, $C_2$ denote two closed embedded curves
with constant geodesic curvature. Assume that the first
eigenvalues for the Jacobi operator~\eqref{eq:jacobi} satisfy
$\lambda_1(C_1)\le0$ and $\lambda_1(C_2)< 0$.
We will see that $C_1\cup C_2$ is an unstable curve.

Fix $\phi_{1}$ (resp. $\phi_{2}$) an eigenfunction
associated to $\lambda_1(C_1)$ (resp. $\lambda_1(C_2)$).
Then we have $J(\phi_{1}) + \lambda_1(C_1)\,\phi_{1}=0$.
Take $\alpha\in\rr$ such that
\[
\alpha\,\int_{C_1}\phi_{1} + \int_{C_2}\phi_{2}=0.
\]
Hence the function
\begin{equation*}
u=\left\{\begin{tabular}{l}
    $\alpha\,\phi_{1}$, \quad in $C_1$,
    \\
    $\phi_{2}$, \quad in $C_2$
    \end{tabular}\right.
\end{equation*}
has mean zero and gives
\[
I(u)=\alpha^2\,\lambda_1(C_1)\,\int_{C_1}\phi_{1}^2 +
\lambda_1(C_2)\,\int_{C_2}\phi_{2}^2 < 0,
\]
so $C_1\cup C_2$ is not stable.

We know that the first eigenvalue for the Jacobi operator
is negative in unduloid type curves,
by Lemma~\ref{le:valorpropio}, and in nodoid type curves,
by Remark~\ref{re:eigennodoid}. Moreover,
for a parallel $\mathbb{S}^1\times\{t\}$,
it is easy to check that
\[
\lambda_1(t)=-\frac{(f')^2 - f f''}{f^2}(t)=-(K+h^2)(t).
\]
Hence $\lambda_1(t)$ is negative if $(K+h^2)(t)>0$,
and vanishes if $(K+h^2)(t)=0$.
Finally, for any vertical geodesic $C$, $f|_C$ is
a positive eigenfunction with zero eigenvalue, and
so $\lambda_1(C)=0$.

In case an unduloid type curve belongs to $\ptl\Om$,
it follows from the above arguments that
necessarily $\Om$ is a set bounded by such a curve and a circle of
revolution $\sph^1\times \{t\}$ with positive first eigenvalue,
equivalently $(K + h^2)(t)<0$, and then $\Om$ is a region of type v).

In case a nodoid type curve belongs to $\ptl\Om$, an analogous
reasoning will give that $\Om$ is a region of type i)
(by Lemma~\ref{le:nodsim}) or of type vi).

Assume now that the boundary of $\Om$ does not contain neither an
unduloid type curve nor a nodoid type curve. By Lemma~\ref{le:av},
any region of type iv) is stable, and if circles of revolution
appear in $\ptl\Om$, the remaining possibilities are regions of
type ii), by Lemma~\ref{le:stabannu}, and of type iii),
by Lemma~\ref{le:curvatura} and the first observation
above in the proof.
\end{proof}

\begin{remark}
We will now check that sets of type $v)$,
annuli bounded by an unduloid type curve
and a circle of revolution,
appear as stable regions in some surfaces.
Consider a standard torus of revolution
of Example~\ref{ex:standard}, with $r<a<3\,r$.
Denote by $\ga_T$ the unduloid type curve
$\ga(\theta,T,h_o(T))$, with initial condition
$T$ close to $\tilde{t}=\pi r/2$.
By Lemma~\ref{le:stableunduloids}, these curves are stable.
Let $\Omega_T$ be the annulus bounded by $\ga_T$,
and the corresponding circle of revolution
$\mathbb{S}^1\times\{t(T)\}$, with
positive first eigenvalue (that is, contained in
$K+h^2<0$) and with the same geodesic
curvature with respect to the inner normal.
Let us prove that $\Om_T$ is stable.

Consider $u$ any mean zero function defined
on $\partial\Omega_T$, normalized so that
$\int_{\ptl\Om} u^2=1$, and let $u_1$, $u_2$ be
the restrictions of $u$ to $\ga_T$ and $\mathbb{S}^1\times\{t(T)\}$,
respectively. It is possible to express $u_i=c_i + v_i$,
with $c_i$ a real constant and $v_i$ a mean zero function,
for $i=1,2$. We have
\begin{align*}
&I(u)=I(u_1) + I(u_2)=
\\
&-c_1^2\,\int_{\ga_T}(K+h^2)\,-2\,c_1\int_{\ga_T}v_1\,(K+h^2) -c_2^2\,L(t(T))\,(K+h^2)(t(T))+
\\
&\quad I(v_1,v_1) + I(v_2,v_2).
\end{align*}

Observe that $I(v_i,v_i)\geq 0$ by the stability of each boundary curve.
Moreover, the mean zero condition on $u$ gives $c_1\,L(\ga_T)=-c_2\,L(t(T))$,
where $L(\ga_T)$ denotes the lenght of $\ga_T$, and so
\begin{align}
\label{Omega}
I(u)\geq &-c_1^2\,L(\ga_T)^2\,\bigg(\frac{\int_{\ga_T} (K+h^2)}{L(\ga_T)^2} +
\frac{(K+h^2)(t(T))}{L(t(T))}\bigg) -
\\
&2\,c_1\,\int_{\ga_T}v_1\,(K+h^2).\notag
\end{align}

When $T$ is close to $\tilde{t}$, it follows that $\ga_T$ is close to
$\mathbb{S}^1\times\{\tilde{t}\}$ and, consequently,
$\Omega_T$ is close to the nonsymmetric annulus associated
to $\tilde{t}$, namely $\mathbb{S}^1\times [t',\tilde{t}]$.
Since this annulus satisfies \emph{strictly}
the stability condition ~\eqref{eq:nosim},
we conclude that the first summand in \eqref{Omega} is positive
for $T$ close enough to $\tilde{t}$.

On the other hand, denoting by $v_1^{+}=\max\{v_1,0\}$ and $v_1^{-}=-\min\{v_1,0\}$
we have that $v_1=v_1^{+}-v_1^{-}$, and so
\begin{align*}
\int_{\ga_T} v_1\,(K+h^2)=&\int_{\ga_T} v_1\,K=\int_{\ga_T} v_1^{+}\,K - \int_{\ga_T} v_1^{-}\,K \leq
\\
&K(s_1)\,\int_{\ga_T} v_1^{+} - K(s_2)\,\int_{\ga_T}v_1^{-},
\end{align*}
where $K(s_1)=\max_{\ga_T}K$ and $K(s_2)=\min_{\ga_T}K$.
Since $\ga_T$ lies in a narrow band around $\mathbb{S}^1\times\{\tilde{t}\}$,
we have $K(s_1)-K(s_2)=\varepsilon_{T}$,
so that $K(s_1)=K(s_2)+\varepsilon_{T}$ and so
\[
\int_{\ga_T} v_1\,(K+h^2)\leq\varepsilon_{T}\,\int_{\ga_T} v_1^{+}.
\]

\begin{figure}[ht]
\centerline{\includegraphics[width=0.6\textwidth]{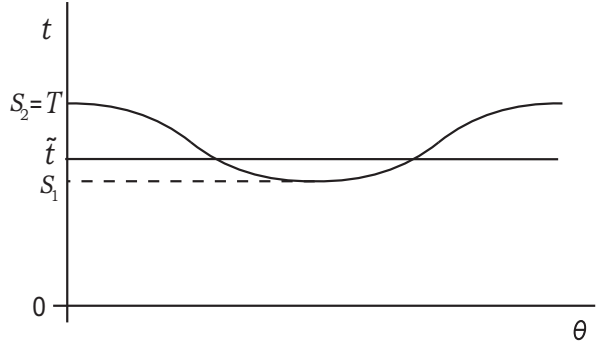}}
\caption{The unduloid type curve $\ga_T$ and
the circle of revolution $\mathbb{S}^1\times \{\tilde{t}\}$}
\label{fig:Omegat}
\end{figure}
Hence
\[
\bigg|\int_{\ga_T} v_1\,(K+h^2)\bigg| \leq \varepsilon_{T}\,\bigg|\int_{\ga_T} v_1^+\bigg| \leq
\varepsilon_{T}\bigg(\int_{\ga_T} (v_1^{+})^2\bigg)^\frac{1}{2}\,L(\ga_T)^\frac{1}{2}.
\]
Since
\[
\int_{\ga_T}(v_1^{+})^2\leq\int_{\ga_T}v_1^2 \leq\int_{\ga_T}u_1^2 = \int_{\ga_T}u^2\leq 1,
\]
we obtain
\[
\bigg|\int_{\ga_T} v_1\,(K+h^2)\bigg|\leq \varepsilon_{T}\,L(\ga_T)^{\frac{1}{2}},
\]
so that the second summand in \eqref{Omega} is negligible,
when $T$ is close enough to $\tilde{t}$. We conclude that $I(u)$
is positive and hence $\Omega_T$ is stable, for $T$ close to $\tilde{t}$.
\end{remark}

\begin{remark}
\label{re:tororaro} We will now show that sets of type vi),
unions of a disk and a symmetric annulus,
actually occur as stable regions in some surfaces.
Let $\widetilde{M}$ be one of the surfaces of the family
described in Example~\ref{ex:tororaro}, obtained from
a sphere and two hyperbolic annuli.

Given a disk $D\subset\widetilde M$, bounded by
a closed embedded nodoid type curve, and a symmetric
annulus $B\subset\widetilde M$ with $h(\ptl D)=h(\ptl B)=h$ with
respect to the inner normal,
following lemma states the stability condition for $\Om=D\cup B$.

\begin{lemma}
\label{le:ejemplo} In the above conditions, $\Om=D\cup B\subset\widetilde{M}$ is
stable if and only if $D$ and $B$ are stable sets, and
\begin{equation}
\label{eq:ejemplo}
\frac{(a^2 + h^2 )^{3/2}}{2\pi} -
\frac{(b^2 - h^2)^{3/2}}{4\pi b c} \leq 0,
\end{equation}
where $a^2$ and $-\,b^2$ are respectively the Gauss curvatures of
the initial sphere and hyperbolic annuli, and $2\pi c$ is the
minimum length of a circle of revolution contained in the
hyperbolic piece.
\end{lemma}

\begin{proof}
Suppose first that $D$ and $B$ are stable and \eqref{eq:ejemplo}
holds. Let $u$ be a mean zero function defined on $\ptl\Om$,
and name $u_1,\,u_2$ the restrictions of $u$ to $\ptl D,\,\ptl B$, respectively.
It is possible to express $u_i=c_i + v_i$, 
with $c_1,\,c_2$ real constants and $v_1,\,v_2$ mean zero
functions on $\ptl D$ and $\ptl B$. These decompositions easily
allows to check that
\begin{align*}
I(u)&=I(u_1)+I(u_2)\geq I(c_1) + I(c_2)
\\
&=-L(\ptl D)\,c_1^2\,(a^2 + h^2) - L(\ptl B)\,c_2^2\,(-b^2+h^2),
\end{align*}
since $v_i$ has mean zero and $D$ and $B$ are stable, and where $L$
denotes the length. The mean zero value condition of $u$ gives
$c_1\,L(\ptl D) = -c_2\,L(\ptl B)$, and then
\[
I(u)\geq -L^2(\ptl D)\,c_1^2\,\bigg(\,\frac{a^2 + h^2}{L(\ptl
D)} + \frac{-b^2 + h^2}{L(\ptl B)}\,\bigg).
\]

By Lemma~\ref{le:stabannu} we know that $B=\mathbb{S}^1\times[-t,t]$
will be contained in the hyperbolic piece since it is stable. Then
\[
L(\ptl B)=2\,L(t)=4\pi c\,\cosh(d-bt)=4 \pi b c\,(b^2-h^2)^{-1/2},
\]
where in the last equality we have used that $h(t)=h$.
Straightforward calculations give that
$L(\ptl D)=2 \pi\,(a^2 + h^2)^{-1/2}$.
Therefore
\begin{equation}
\label{fin} I(u)\geq -L^2(\ptl D)\,c_1^2\,\bigg(\,\frac{(a^2 +
h^2)^{3/2}}{2\,\pi} - \frac{(b^2 -
h^2)^{3/2}}{4 \pi b c}\,\bigg).
\end{equation}
As we are assuming that \eqref{eq:ejemplo} is verified, then
\eqref{fin} is nonnegative, and so $\Om$ is stable.

Assume now that $\Om$ is a stable region. Then clearly $D$ and $B$
are also stable. Consider the mean zero function
\begin{equation}
\label{eq:const} u=\left\{\begin{tabular}{l}
    $L(\ptl D),\quad\ptl B,$
    \\
    $-L(\ptl B),\quad\ptl D.$
    \end{tabular}\right.
\end{equation}
Using this function in the index
form, we will have $I(u)\geq 0$ by the stability of $\Om$, which
trivially gives \eqref{eq:ejemplo}.
\end{proof}

The above lemma gives us the existence of surfaces where
sets of type vi) are stable. For instance, consider
the surface $\widetilde M$ given by $a=1$, $t^*=\pi/6$,
with $b=0.578$ and $c=0.0410512$. Then it can be checked
that, for geodesic curvature $h=0.4$, there exists a stable
union of a disk and a symmetric annulus.
\end{remark}

\section{Isoperimetric regions}
\label{sec:isoperimetric}

Since we are studying compact surfaces, well-known
results from geometric measure theory \cite{m} ensure the
existence of isoperimetric solutions for any value of the area.
Moreover, as any isoperimetric
region is stable, the candidates are given by
Theorem~\ref{th:main}. We first discard sets of tipe iv)
composed by several vertical annuli.

\begin{lemma}
The union of two or more vertical annuli is not an isoperimetric
region.
\end{lemma}

\begin{proof}
Simply rotate a vertical annulus until meeting another one and,
after eliminating the duplicated vertical geodesic, we get a
less-perimeter region enclosing the same area.
\end{proof}

The following lemma shows that when the boundary of an
isoperimetric region is an unduloid type curve and a circle of
revolution, some restrictions appear. Call \emph{lower half}
of $M$ to $\sph^1 \times [-t_0,0]$ and \emph{upper half} to
$\sph^1 \times [0, t_0]$.

\begin{lemma}
\label{le:restricciones}

Let $\Om$ be an isoperimetric region in the warped product $\sph^1
\times [-t_0,t_0]$ bounded by an unduloid type curve $C$ and
a circle of revolution $\sph^1\times \{t\}$. Then
\begin{itemize}
\item[i)] $C$ and $\sph^1 \times \{t\}$ are not contained
in the same half.

\item[ii)] $C$ does not intersect $\sph^1 \times \{-t\}$.
\end{itemize}
\end{lemma}

\begin{proof}
Without loss of generality we can assume that $t\in [-t_0,\, 0]$.
Since $\Om$ has to be stable, by Theorem~\ref{th:main} we have
$(K+h^2)(t)< 0$. Consequently, $t\in [-t_0,t_c]$.

i) Suppose both curves lie in the half $\sph^1 \times
[-t_0,\, 0]$. Call $t_m,\,t_M$ the minimun and the maximum of
$t|_C$, respectively. If $t = -t_0$, it can be checked by
using \eqref{eq:primeraintegral} that $f(t_m)=f(t_M)$, and then
$C$ intersects symmetrically $\sph^1\times\{0\}$. This
contradicts the fact that $C$ is contained in the lower
half. Hence $t \neq -t_0$. Since $C$ has to intersect
$(f')^2 - f f'' = 1$, it will be necessarily contained in $\sph^1
\times [t,\, 0]$. But then we can construct a new set
enclosing the same area with \emph{strictly} less perimeter: consider a parallel
$\sph^1 \times \{t^*\}$ intersecting $C$, with $t^*$ close
to $t_M$, and replace the piece of $C$ above the parallel by
the corresponding segment of parallel.
Also replace $\sph \times \{t\}$ by $\sph \times
\{\overline t\}$ with $ \overline t < t$. There exist appropriate
$t^*$ and $\overline t$ for which the new set encloses the same
area and it can be checked that it has less perimeter, which is
contradictory since $\Om$ is isoperimetric.

ii) Suppose $C$ intersects $\sph^1 \times \{-t\}$. Replace the
piece of $C$ above the parallel $\sph^1 \times \{-t\}$ by the
corresponding segment of the parallel. By reflecting the replaced
piece with respect to $\sph^1\times\{0\}$ and removing a segment
of $\sph^1 \times \{t\}$, we obtain a new set with the same
perimeter and enclosing the same area that $\Om$, but with no
regular boundary, which is a contradiction.
\end{proof}

\begin{remark}
Assume that a disk with constant geodesic boundary,
and contained in a region with constant Gauss curvature,
is a component of an isoperimetric region.
Then, by using a comparison argument (\cite{ritore}, Lemma~2.7),
it can be checked that such a constant is precisely
the maximum of the Gauss curvature.
\end{remark}

\subsection{Isoperimetric regions in standard tori}
We will now focus on standard tori described in Example~\ref{ex:standard}.
Denote by $\beta$ the total area of a given torus.
Fixing the value of $a$, we have studied
the isoperimetric solutions for the different values of $r$.

\begin{theorem}
\label{th:standard} Consider a standard torus,
parametrized by $r$ and $a$, with $r\leq a/2$.
Then the isoperimetric regions are
disks with constant geodesic curvature, and vertical annuli,
and their complements.
\end{theorem}

\begin{proof}
Since $r\leq a/2$, it is easy to see that the
length of any circle of revolution, and
therefore of any unduloid type curve, is
greater than or equal to the length of a
vertical geodesic. This fact reduces the
candidates to disks and vertical annuli,
which actually appear as isoperimetric.
\end{proof}

When $r>a/2$ we have observed, with numerical computations,
different behaviors as $r$ increases.
For values of $r$ close to $a/2$,
it happens as in Theorem~\ref{th:standard}; there are also values of
$r$ for which disks (or complements of them) are
solutions for any quantity of the area; then, there exists an
interval for $r$ where the isoperimetric regions are
disks, symmetric annuli afterwards, and finally disks
for values of the area near $\beta/2$; and for $r$ close enough to
$a$, we get disks, symmetric annuli, nonsymmetric annuli and
disks again as solutions.

\begin{figure}[ht]
\centering{
\subfigure[$r=0.77\,a$]{\label{perfil0}\includegraphics[width=0.44\textwidth]{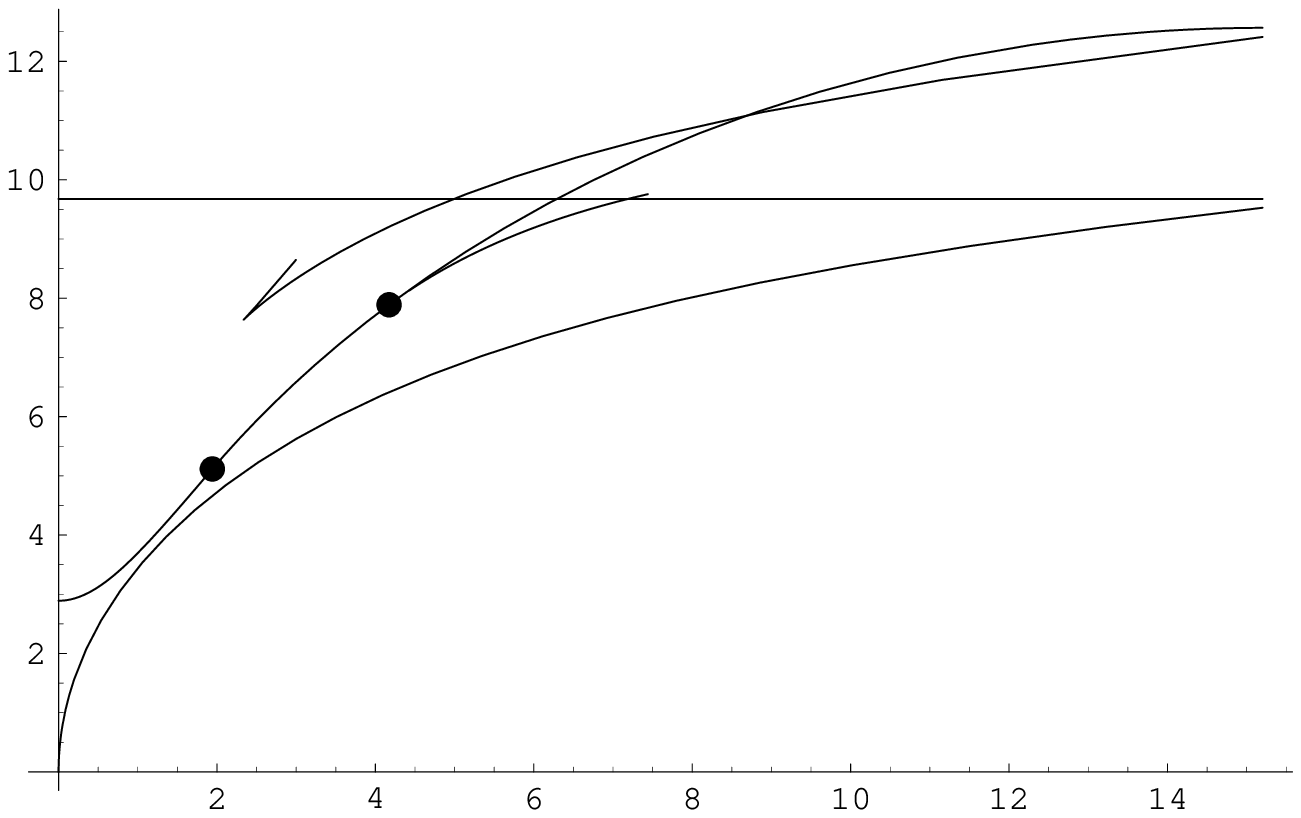}}
\hspace{0.1\textwidth}
\subfigure[$r=0.9\,a$]{\label{perfil2}\includegraphics[width=0.44\textwidth]{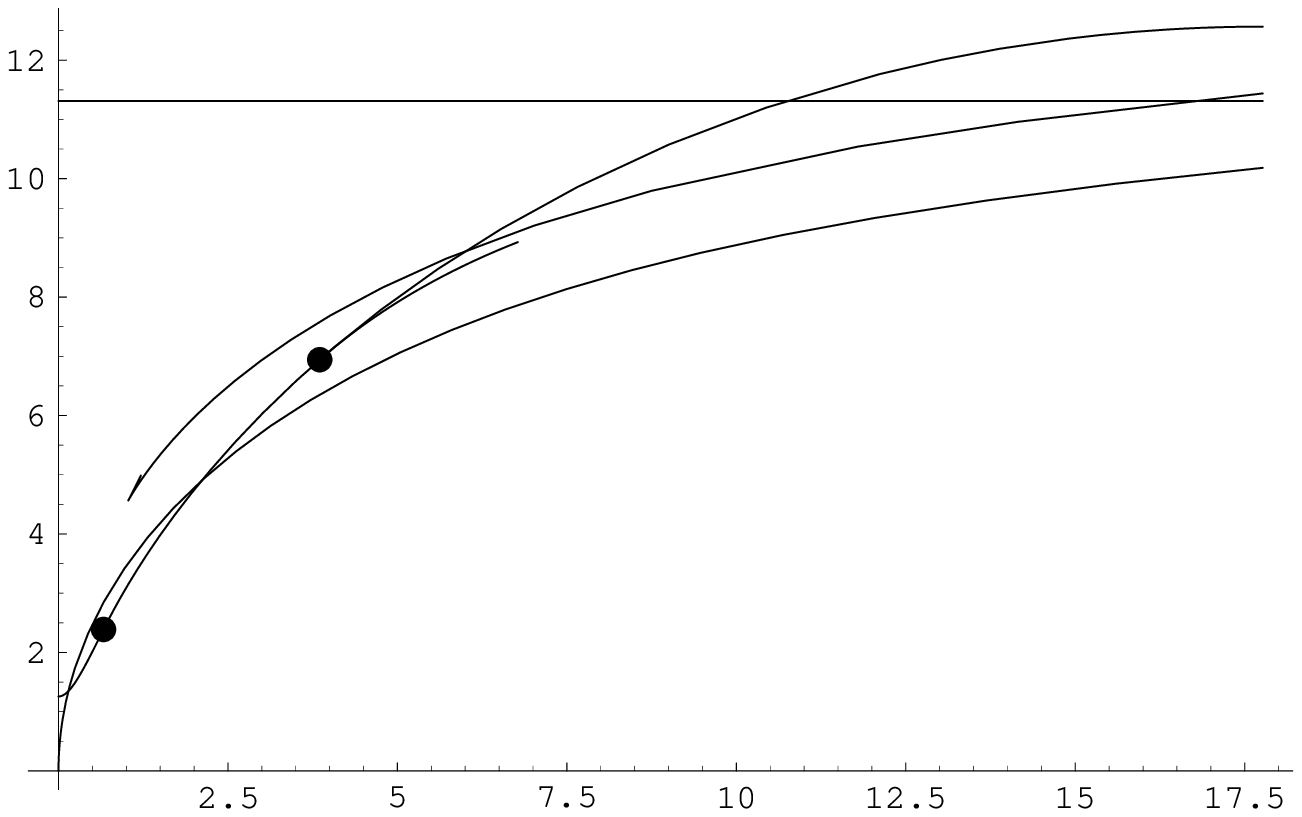}}
}
\caption{Graphs showing the perimeter of each candidate in two different
standard tori}
\label{fig:grafos}
\end{figure}

The perimeter of the isoperimetric candidates,
as a function of the area enclosed,
is depicted in two different standard tori of revolution
in Figure~\ref{fig:grafos}.
In each graph, the curve starting from the origin
corresponds to the perimeter of disks,
and the above one not touching the vertical axis
corresponds to the perimeter of unions of a disk and
a symemtric annulus.
In the first graph, we observe that the perimeter of disks
is less than any other set, for any quantity of area.
However, in the second one, symmetric and nonsymmetric annuli
are also solutions (the black dots show the transition between symmetric and
nonsymetric annuli, and between nonsymmetric annuli and sets of type v)).

\begin{remark}
\label{re:no}
We have not found any stardard torus where a domain bounded by
an unduloid type curve and a circle of revolution is an isoperimetric
region.
\end{remark}

\begin{remark}
\label{re:grandesdiscos} We have observed numerically that
for any standard torus, the isoperimetric region
for area $\beta/2$ is either a disk with constant geodesic
boundary, or a vertical annulus.
We will see that this fact does not hold in others surfaces.
\end{remark}

\subsection{Isoperimetric regions in $\widetilde M$}
Consider now the surfaces $\widetilde{M}$ introduced in
Example~\ref{ex:tororaro} and call as above $\beta$ the total
area of the surface. We have studied carefully the isoperimetric
problem in this case, specially when vertical annuli have too much
perimeter (that is, when hyperbolic annuli are long enough), and hence
they are not isoperimetric solutions.
In this setting, we have seen that isoperimetric regions
are disks bounded by constant geodesic curvature curves
for small areas (entirely contained in the spherical piece,
where the Gauss curvature achieves constantly its maximum),
and afterwards symmetric annuli (with its boundary contained
in the hyperbolic piece); or disks, all symmetric annuli which are stable, and finally
nonsymmetric annuli for quantities of area near $\beta/2$ (with a
circle of revolution of its boundary in the spherical piece, and
another one in the hyperbolic part).

Annuli bounded by an unduloid type curve and
a circle of revolution also appear as
isoperimetric regions in these surfaces.
Consider a nonsymmetric annulus which is
an isoperimetric region in $\widetilde{M}$.
By rotating slightly the spherical piece,
keeping unchanged the hyperbolic annuli,
we will obtain another isoperimetric region,
since area and perimeter are preserved.
This new region is now bounded by an unduloid
type curve (contained in the spherical part)
and a circle of revolution (in the hyperbolic one).

Then we have
\begin{theorem}
\label{th:tororaro}
Let $\widetilde{M}$ be one of the surfaces
described in Example~\ref{ex:tororaro}, obtained from a
sphere and a hyperbolic annulus.
Then the isoperimetric regions may be disks
with constant geodesic boundary, symmetric
or nonsymmetric annuli, annuli bounded by an unduloid
type curve and a circle of revolution, unions of
a disk and a symmetric annulus, vertical annuli
or the complement of one of these sets.
\end{theorem}

\begin{remark}
In \cite[Prop.~3.4.]{pr} it is proved that in high dimensions,
there are isoperimetric domains bounded by hypersurfaces of
revolution generated by unduloid type curves. Up to our knowledge,
this kind of curves had not appear as part of the boundary of an
isoperimetric solution in any surface. The above result shows that
this fact may occur.
\end{remark}

\begin{remark}
Note that for a surface $\widetilde{M}$ of Example~\ref{ex:tororaro},
symmetric annuli, nonsymmetric annuli and annuli bounded by
an unduloid type curve and a circle of revolution
may be isoperimetric solutions for area $\beta/2$.
\end{remark}

\begin{remark}
Although a set consisting of a disk and a symmetric
annulus can be stable (Lemma~\ref{le:ejemplo}),
we have not found any surface where such
a region is isoperimetric. For instance, this kind of sets does
not appear in most of the surfaces $\widetilde M$,
since the constant geodesic curvature condition is hardly
verified. However, it seems that only a length comparison argument
might discard them as isoperimetric.
\end{remark}


\begin{thebibliography}{\bf HHHH}

\bibitem[{\bf BB}]{bb}
L. Barbosa, P. B\'erard, \emph{Eigenvalue and ``twisted" eigenvalue
problems, applications to CMC surfaces}, J. Math. Pures Appl.~\textbf{9}
(2000), no.~5, 427--450. \MR {2001f:58064}

\bibitem[{\bf BGS}]{bgs}
J. L. M. Barbosa, J. M. Gomes and A. M. Silveira,
\emph{Foliation of 3-dimensional space forms by
surfaces with constant mean curvature}, Bol. Soc. Brasil. Mat.~\textbf{18}
(1987), no.~2, 1--12. \MR {90j:53054}

\bibitem[{\bf BP}]{bapa}
C. Bavard and P. Pansu, \emph{Sur le volume minimal de {${\bf
R}\sp 2$}},  Ann. Sci. {\'E}cole Norm. Sup. (4)~\textbf{19}
(1986), no.~4, 479--490. \MR {88b:53048}

\bibitem[{\bf BC}]{bc}
I. Benjamini and J. Cao, \emph{A new isoperimetric comparison
theorem for surfaces of variable curvature}, Duke Math. J.
\textbf{85} (1996), 359--396. \MR {97m:58046}

\bibitem[{\bf Ch}]{chavel}
I. Chavel, \emph{Eigenvalues in Riemannian Geometry},
Academic Press, 1984. \MR {86g:58140}

\bibitem[{\bf CL}]{cl}
E.A. Coddington and N. Levinson, \emph{Theory of Ordinary Differential Equations},
McGraw-Hill, New York, 1955. \MR {16,1022b}

\bibitem[{\bf CHLL}]{dp}
J. Corneli, P. Holt, G. Lee, N. Leger, E. Schoenfeld and B.
Steinhurst, \emph{The double bubble problem on the flat
two-torus}, Trans. Amer. Math. Soc. \textbf{356} (2004), no.~9,
3769--3820. \MR {2005b:53011}

\bibitem[{\bf H}]{H}
H. Howards, \emph{Soap bubbles on surfaces}, undergraduate thesis,
Williams College, 1992.

\bibitem[{\bf HHM1}]{hhm1}
H. Howards, M. Hutchings and F. Morgan, \emph{The isoperimetric
problem on surfaces}, Amer. Math. Monthly \textbf{106} (1999), no.~5,
430--439. \MR {2000i:52027}

\bibitem[{\bf HHM2}]{hhm2}
H. Howards, M. Hutchings and F. Morgan, \emph{The isoperimetric
problem on surfaces of revolution of decreasing Gauss curvature},
Trans. Amer. Math. Soc. \textbf{352} (2000), no.~11, 4889--4909.
\MR {2001b:58024}

\bibitem[{\bf HL}]{hl}
W.~H. Huang, C.~C. Lin, \emph{Negatively curved sets on surfaces of
constant mean curvature in $\rr^3$ are large}, Arch. Rational Mech. Anal.
\textbf{141} (1998), no.~2, 105--116. \MR {99j:58049}

\bibitem[{\bf K}]{koiso}
M. Koiso, \emph{Deformation and stability of surfaces with
constant mean curvature}, Tohoku Math. J. \textbf{54} (2002), 145--159.
\MR {2003j:58021}

\bibitem[{\bf M}]{m}
F. Morgan, \emph{Geometric Measure Theory: a Beginner's Guide},
third ed., Academic Press Inc, 2000. \MR {2001j:49001}

\bibitem[{\bf Mi}]{mi}
J. Milnor, \emph{Morse theory}, Princeton Univ. Press, Princeton,
N.J., 1963. \MR {29:634}

\bibitem[{\bf O}]{O}
R.~Osserman, \emph{Bonnesen-style isoperimetric inequalities},
Amer. Math. Monthly \textbf{86} (1979), no.~1, 1--29.
\MR {80h:52013}

\bibitem[{\bf P}]{p}
P.~Pansu, \emph{Sur la r\'egularit\'e du profil
isop\'erim\'etrique des surfaces riemanniennes compactes}, Ann.
Inst. Fourier (Grenoble) \textbf{48} (1998), no.~1, 247--264. \MR
{99i:53035}

\bibitem[{\bf PR}]{pr}
R. H. L. Pedrosa, M. Ritor{\'e}, \emph{Isoperimetric domains in
the Riemannian product of a circle with a simply connected space
form and applications to free boundary problems}, Indiana Univ.
Math. J. \textbf{48} (1999), no.~4, 1357--1394. \MR {2001k:53120}

\bibitem[{\bf R}]{ritore}
M. Ritor{\'e}, \emph{Constant geodesic curvature curves and
isoperimetric domains in rotationally symmetric surfaces}, Comm.
Anal. Geom. \textbf{9} (2001), no.~5, 1093--1138. \MR
{2003a:53018}

\bibitem[{\bf S}]{s}
E. Schmidt, \emph{\"{U}ber eine neue Methode zur Behandlung einer
Klasse isoperimetrischer Aufgaben im Grossen}, Math. Z.
\textbf{47} (1942), 489--642. \MR{7} (1946)

\bibitem[{\bf T}]{t}
P.~Topping, \emph{The isoperimetric inequality on a surface},
Manuscripta Math. \textbf{100} (1999), no.~1, 23--33.
\MR{2000i:53107}

\end{thebibliography}

\providecommand{\bysame}{\leavevmode\hbox
to3em{\hrulefill}\thinspace}
\providecommand{\MR}{\relax\ifhmode\unskip\space\fi MR }
\providecommand{\MRhref}[2]{%
   \href{http://www.ams.org/mathscinet-getitem?mr=#1}{#2}
} \providecommand{\href}[2]{#2}

\end{document}